\documentclass[preprint,showpacs,preprintnumbers,nofootinbib,amsmath,amssymb,aps,prc,longbibliography]{revtex4-1}
\usepackage{graphicx}
\usepackage{dcolumn}
\usepackage{bm}
\begin{document}
\title{Standard canvas and stretcher sizes satisfying golden and silver ratios as well as optimal use of material}
\author{Nguyen Dinh Dang}
  \email{Email: ndinhdang@gmail.com, URL: http://ribf.riken.go.jp/~dang/}

\date{\today}
\begin{abstract}
The sizes of canvases and stretchers for oil painting have been standardized in France in 19th century and widely accepted in many countries in Europe as well as in Japan so far. These standard sizes do not follow the golden ratio and porte d'harmonie number as has been often claimed. In this work, a general formula is derived to calculate the sizes of painting canvases and stretchers, which satisfy exactly the golden ratio and the porte d'harmonie number as well as the rule for optimal use of material in their mass production. Based on this formula new canvas and stretcher sizes are calculated for all three figure, landscape and marine formats.
\end{abstract}

                              
\maketitle
\section{Introduction}
\label{Intro}
Canvas has been used in oil painting starting from the end of 15th century in Italy. A century later the use of canvas in painting spread to England, where canvases and stretchers were standardized and manufactured to be sold to artists~\cite{England}. In 19th century a new system of standards for canvas sizes was introduced in France (the so-called second standardization after the first one in 17th century)~\cite{Labreuche1,Labreuche2}. This system has been accepted throughout Europe and even imported to Japan, where it remains the only one for canvas and strainer formats available from the local manufacturers of art materials. In the French standards the size increment goes discretely from size 0 (called point 0 or toile de 0) to size 120 (called point 120 or toile de 120) with 20 sizes for each of three formats, which are the figure (F), landscape (paysage, P) and marine (M) ones, as shown in Table \ref{FrenchTable}. The Japanese standard sizes are identical to the French ones except some slight errors due to conversion and rounding up from the metric units to the Japanese shaku units in the Meiji era and the reverse conversion back to the metric units in 20th century (1 shaku $=$ 10/33 m $\simeq$ 30.3 cm). Beyond that, the Japanese standards have also been extended to sizes 130, 150, 200, 300, and 500. 
\begin{table}
\begin{center} 
    \caption
        {French standard canvas sizes and their Japanese extension for sizes 130 -- 500, which have been corrected for errors due to conversion between metric and shaku units. The heights for P130 and M130, absent in the Japanese extension, have also been added by using Eq. (\ref{optimal}).\label{FrenchTable}}
       \vspace{2mm}        
\begin{tabular}{|c|c|c|ccc|}
\hline\hline 
~~~~{No}~~~~ &~~~~{Point}~~~~ & ~~~{Width}~~~ &\multicolumn{3}{c|}{~~~~~~Height (cm)~~~~~~}\\
&&(cm)&\multicolumn{1}{c}{~~~~F~~~~}&\multicolumn{1}{c}{~~~~P~~~~}&\multicolumn{1}{c|}{~~~~M~~~~}\\
\hline
1 &0 &       18&   14&	12&	10\\
2 &1 &	22&   16&	        14&	12\\
3 &2 &	24&   19&	         16&	14\\
4 &3 &	27&    22&        19&	16\\
5 &4 &	33&   24&	         22&	19\\
6 &5 &	35&    27&	24&	22\\
7 &6 &	41&    33&	27&	24\\
8 &8 &	46&    38&	33&	27\\
9 &10&	55&    46&	38&	33\\
10 &12&	61&    50&	46&	38\\
11 &15&	65&    54&	50&	46\\
12 &20&	73&    60&	54&	50\\
13 &25&	81&    65&	60&	54\\
14 &30&	92&    73&	65&	60\\
15 &40&	100&  81&	73&	65\\
16 &50&	116&  89&	81&	73\\
17 &60&	130&  97&	89&	81\\
18 &80&	146&  114	&      97&	89\\
19 &100&	162&  130	&    114&	97\\
20 &120&	195&  130	&    114&	97\\
21 &130& 195&  162&    130& 114\\
22 &150&  228&  182&   162&  146\\
23 &200&  260&  195&    182&  162\\
24 &300& 292&   219&    197&    182\\
25 &500& 334&   250&    219&    197\\
 \hline\hline
\end{tabular} 
\end{center}
\end{table}

Opinions regarding the French canvas standard sizes have been divided. A number of the so-called {\it les initi\'{e}s} insist that these standards haven been derived from the golden ratio $\varphi =(1+\sqrt{5})/2=$ 1.6180339887..., and $\sqrt{2}$, which is also called the number of {\it porte d'harmonie} based on the name given by the French painter Paul S\'{e}rusier (1864 - 1927). Since the silver ratio is given as $\delta_{S} = 1+\sqrt{2}$, one may say that $\sqrt{2}$ is related to the silver number as $\delta_S-1$~\cite{dynamic}. Just it has been claimed that in the standard French canvas sizes, the height (the shorter length) of a canvas in the F format is equal to its width (the longer length) divided by 2$\Phi$ or 2$/\varphi$, where $\Phi$ is  golden ratio conjugate, $\Phi=(\sqrt{5}-1)/2=1/\varphi=$ 0.6180339887.... As for the P format, its height is equal to its width divided by $\sqrt{2}$. Finally, the height of M format  is equal to its width divided by $\varphi$. This allows two identical F format canvases to be jointed along their widths to make one M format canvas, whose width is twice the height of the F format canvas. Just to quote some references, this interpretation is given on the website of a prestigious art university~\cite{MAU}, and on the technical information page of one of the world's largest manufacturers of professional art materials~\cite{Holbein}. 

\begin{figure}
       \includegraphics[width=9cm]{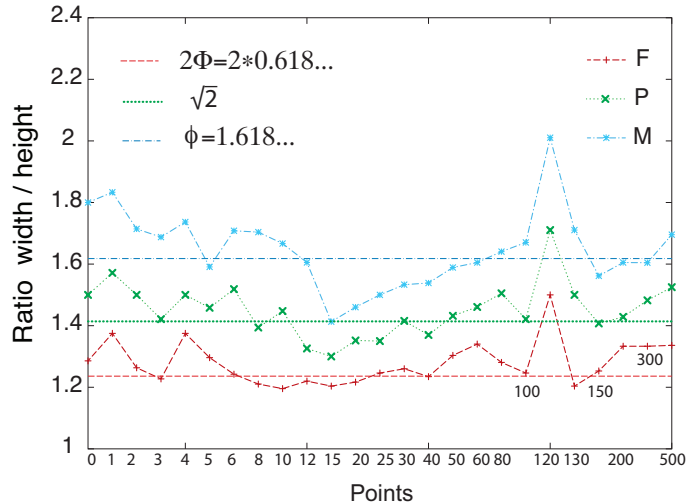}
       \caption{(Color online)  Width-to-height ratios as functions of points in French standard canvas formats.  The (red) crosses, (green) saltires, and
       (blue) asterisks denote the ratios corresponding to the  F, P, and M formats, respectively. The lines connecting these symbols are drawn to guide the eye. The horizontal straight lines stand for 2$\Phi=2/\varphi$ (red dashed), $\sqrt{2}$ (green dotted), and $\varphi$ (dash-dotted). \label{ratioF}}
    \end{figure}

However, these claims are far from reality. The actual sizes of French standard canvas formats including their Japanese extension, as shown in the Table \ref{FrenchTable}, in general do not satisfy the golden and silver ratios, as can be seen in Fig. \ref{ratioF}. As a matter of fact, from 25 sizes in each of F, P, and M formats, ranging from size 0 to size 500, only 7 sizes (less than 30\%) have the ratios of width to height that match the golden or silver ratio. Indeed, for the F format, only the sizes F3, F6, F12, F25, F40, F100, and F150 have the width-to-height ratios close to 2$\Phi$. For the P format, the sizes P3, P8, P30, P50, P100, P150, and P200 have these ratios close to $\sqrt{2}$, whereas for the M format, the sizes M5, M12, M50, M60, M80, M200, and M300 are those whose ratios of width to height are close to $\varphi$. For the remaining sizes these ratios deviate from the corresponding golden and silver ratios, and the deviations are particularly large in some cases like at sizes 1, 120 and M15.

The first original and systematic research on the history of French standard canvas formats as well as their subsequent use and development in France has been carried out by Anthea Callen. In her PhD thesis in 1980~\cite{Callen1}, summarized in Ref. \cite{Callen2}, she has demonstrated how these formats did not conform to any harmonious ratios. Her detailed diagrams and results of calculations have revealed that these standard canvas formats were actually the result of economic determinants, such as the optimal use of fabric width produced by the fixed loom size, rather than aesthetic ones. This research has been taken up later on by Pascal Labreuche~\cite{Labreuche1,Labreuche2}, Ann Hoenigswald~\cite{Ann} and others. 

As a matter of fact, as Callen has already observed~\cite{Callen1,Callen2}, the canvas widths of the sizes from 0 to 120 in the French standards are the same as those mentioned by Antoine-Joseph Pernety already in 1757 in Ref.~\cite{Pernety}, where he pointed out that the points assigned to the different sizes of the prepared and pre-stretched canvases just referred to their costs in sols (or sous). For example, a size 20 canvas would cost 20 sous. Although Pernety presented the measurements in the ancient units, pieds and pouces, the later conversion into the metric units (in cm and mm) does not change these sizes significantly except for some rounding up, which remain as the sizes of F format canvases nowadays. In Ref. \cite{Labreuche2} Labreuche has also noted that attempts have failed to find the golden ratio or a close value in the table of standard canvas formats in 18th century.

This point of view seems to hold also because of rationalized wood cutting in manufacturing canvas stretcher bars~\cite{Labreuche2}. Indeed, Table \ref{FrenchTable} reveals a striking rule regarding the heights of the F, P, and M formats in the successive sizes, namely the height of the  $i$th size in the F format is equal to that of the $(i+1)$th size in the P format, which in turn is the same as that of the $(i+2)$th size in the M format. For instance, the heights of the F0, P1, M2 are the same and equal to 14 cm; so are those of the F1, P2, M3, which are equal to 16 cm, etc. This rule allows the manufacturer to make the shorter stretcher bars of the same length for assembling stretchers of 3 formats (F, P, and M) in the same way as the longer stretcher bars of the same length can be used for all three formats. Consequently it reflects the optimal use of material in manufacturing stretcher bars and pre-stretched canvases. In the rest of the present work this rule will be referred to as the rule for optimal use of material. It can be written as
\begin{equation}
F(i) = P(i+1) = M(i+2)~,
\label{optimal}
\end{equation}
where $F(i)$, $P(i+1)$ and $M(i+2)$ denote the heights of the F, P, and M formats, respectively, at their corresponding $i$, $i+1$ and $i+2$ sizes.  An exception is the last size $i=$ 20 (point 120), where the heights are set to be the same as the corresponding heights at $i=19$, because there is no larger size.

Of course, it is trivial to design a new system of standard canvas and stretcher sizes that strictly follow the golden and silver ratios, keeping the same traditional lengths for the widths $W(i)$ at each size $i$. But these new heights, which are obtained by simply dividing the standard widths $W(i)$ by 2$\Phi$, $\sqrt{2}$ and $\phi$ to obtain $F(i)$, $P(i)$, and $M(i)$, respectively, obviously violate the rule for optimal use of material (\ref{optimal}). Golden rectangle canvases with the width-to-height ratio approximating $\varphi$ are available on the market~\cite{market}.

The goal of the present work is to derive new standard canvas and stretcher sizes for canvas, which simultaneously satisfy the golden and silver ratios as well as the rule for optimal use of material.
\section{New standard sizes for canvas and stretchers}
\label{formula}
As has been mentioned in the Introduction, let's index the canvases withs the numbers $i=$ 1, ..., 25 in the order of increasing sizes and call the width, and heights of F, P, and M formats at $i$th size as $W(i)$, $F(i)$, $P(i)$, and $M(i)$, respectively. The table of new standard sizes then has the form of Table \ref{NewTable}.
\begin{table}
\begin{center} 
    \caption
        {The schematic table to be designed. \label{NewTable}}
       \vspace{2mm}        
\begin{tabular}{ccccc}
\hline
~~~~{No}~~~~ & ~~~{Width}~~~ &\multicolumn{3}{c}{~~~~~~Height~~~~~~}\\
&&\multicolumn{1}{c}{~~~~F~~~~}&\multicolumn{1}{c}{~~~~P~~~~}&\multicolumn{1}{c}{~~~~M~~~~}\\
\hline
1 &       W(1)&   F(1) &	P(1) &	M(1)\\
2 &	W(2) &   F(2) &	        P(2)&	M(2)\\
3 &	W(3) &   F(3) &	         P(3) &	M(3)\\
...&	&    &        &	\\
i & W(i) &   F(i) &M(i)&P(i)\\
i+1 & W(i+1) &   F(i+1) &M(i+1)&P(i+1)\\
... &	&   &	&\\
 \hline
\end{tabular} 
\end{center}
\end{table}

These widths and heights should follow the golden ratio and the porte d'harmonie number, as
\begin{equation}
W(i) = 2 F(i)/\varphi~,
\label{golden1} 
\end{equation}

\begin{equation}
W(i) = \sqrt{2} P(i)~,
\label{silver} 
\end{equation}

\begin{equation}
W(i) = \varphi M(i)~,
\label{golden2} 
\end{equation}
as well as the rule for optimal use of material (\ref{optimal}).

From Eqs. (\ref{optimal}) and (\ref{silver}), one obtains
\begin{equation}
W(i+1) =\sqrt{2} P(i+1) = \sqrt{2} F(i)~.
\label{Wi+1}
\end{equation}
Similarly, from Eqs. (\ref{optimal}) and (\ref{golden2}), one obtains
\begin{equation}
W(i+2) =\varphi M(i+2) = \varphi F(i)~.
\label{Wi+2}
\end{equation}
Dividing $W(i+1)$ from Eq. (\ref{Wi+1}) by $W(i)$ from Eq. (\ref{golden1}), one gets the ratio between two widths as
\begin{equation}
\frac{W(i+1)}{W(i)} =\frac{\varphi}{\sqrt{2}}~.
\label{ratio1}
\end{equation}
In the same way, dividing $W(i+2)$ from Eq. (\ref{Wi+2}) by $W(i)$, one obtains
\begin{equation}
\frac{W(i+2)}{W(i)} =\frac{\varphi^2}{2}~.
\label{ratio2}
\end{equation}
Setting the width $W(1)$ of the smallest size ($i=$ 1) equal to $w$, from Eq. (\ref{ratio1}) one has
\begin{equation}
W(2)=W(1)\frac{\varphi}{\sqrt{2}}= w\frac{\varphi}{\sqrt{2}}~,
\label{W2}
\end{equation}
\begin{equation}
W(3)=W(2)\frac{\varphi}{\sqrt{2}}=w\bigg(\frac{\varphi}{\sqrt{2}}\bigg)^2~,
\label{W3}
\end{equation}
\begin{equation}
W(4)=W(3)\frac{\varphi}{\sqrt{2}}=w\bigg(\frac{\varphi}{\sqrt{2}}\bigg)^3~,
\label{W4}
\end{equation}
\[
.......~~~~~~~~~~~~~~~~~~~~~
\]
which also satisfy Eq. (\ref{ratio2}). Proceeding in the same manner, one ends up with
\begin{equation}
W(i)=w\bigg(\frac{\varphi}{\sqrt{2}}\bigg)^{i-1}~.
\label{Wi}
\end{equation}

{\it Proof of Eq. (\ref{Wi})}:

One can prove this formula by induction, namely if Eq. (\ref{Wi}) holds for $i=n$, one has to show that it also holds for $i=n+1$, that is
if
\begin{equation}
W(n)=w\bigg(\frac{\varphi}{\sqrt{2}}\bigg)^{n-1}~,
\label{Wn}
\end{equation}
then
\begin{equation}
W(n+1)=w\bigg(\frac{\varphi}{\sqrt{2}}\bigg)^{n}~.
\label{Wn+1}
\end{equation}
Replacing $W(i)$ in Eq. (\ref{golden1}) with the right-hand side of Eq. (\ref{Wn}), one obtains
 \begin{equation}
F(n)=W(n)\frac{\varphi}{2} = w\frac{\varphi^n}{2(\sqrt{2})^{n-1}}~.
\label{Fn}
\end{equation}
From Eq. (\ref{optimal}) one has $F(n) = P(n+1)$, therefore
 \begin{equation}
P(n+1)= w\frac{\varphi^n}{2(\sqrt{2})^{n-1}}~.
\label{Pn+1}
\end{equation}
 Replacing $P(i=n+1)$ at the right-hand side of Eq. (\ref{silver}) with its expression determined by 
 the right-hand side of Eq. (\ref{Pn+1}), one gets
\begin{equation}
W(n+1)= \sqrt{2}P(n+1) = w\bigg(\frac{\varphi}{\sqrt{2}}\bigg)^n~,
\label{proof}
\end{equation}
which is nothing but Eq. (\ref{Wn+1}) that one has to prove.

Equation (\ref{Wi}) is the central result of this work. It states that the width of the $i$th canvas is equal to that
of the first (smallest) canvas multiplied by the ratio of the golden number ($\varphi$) to $\sqrt{2}$ (the porte d'harmonie number), raised to the $(i-1)$th power. Knowing the width $W(i)$, one can easily calculate the heights $F(i)$, $P(i)$, and $M(i)$ at the same width $W(i)$, which satisfy the golden ratio and the porte d'harmonie number as well as the rule for optimal use of material, by using Eqs. (\ref{golden1}) -- (\ref{golden2}).

Shown in Tables \ref{Table1} and \ref{Table2} are the new standard sizes for canvases and stretchers, calculated from Eqs. (\ref{Wi}) and (\ref{golden1}) -- (\ref{golden2}), where the number in Table \ref{Table2} have been rounded up to mm. They are obtained by using the width for the smallest size. $W(1)=w=$ 18 cm, that is the width of size 0 from the French standard sizes in Table \ref{FrenchTable}.
\begin{table}
\begin{center} 
    \caption
        {New standard sizes for canvas and stretchers. The size names in points in the 2nd column are assigned to the width lengths, which are close to those at the corresponding points in the French standard sizes and their Japanese extension, mentioned in the text. \label{Table1}}
       \vspace{2mm}        
\begin{tabular}{|c|c|c|ccc|}
\hline\hline 
~~~{No}~~~ &~~{Point}~~ & ~~~{Width}~~~ &\multicolumn{3}{c|}{~~~~~~Height (cm)~~~~~~}\\
&&(cm)&\multicolumn{1}{c}{~~~~~~F~~~~~~}&\multicolumn{1}{c}{~~~~~~P~~~~~~}&\multicolumn{1}{c|}{~~~~~~M~~~~~~}\\
\hline
  1 &  0 &   18.000000000 &   14.562306046 &   12.727922279 &   11.124611685\\
  2 &  1 &   20.594211063 &   16.661066905 &   14.562306545 &   12.727922279\\
  3 &  2 &   23.562307183 &   19.062307131 &   16.661067475 &   14.562306545\\
  4 &  3 &   26.958173737 &   21.809620912 &   19.062307784 &   16.661067475\\
  5 &  4 &   30.843462211 &   24.952885347 &   21.809621658 &   19.062307784\\
  6 &  5 &   35.288709483 &   28.549165971 &   24.952886201 &   21.809621658\\
  7 &  6 &   40.374618401 &   32.663752759 &   28.549166948 &   24.952886201\\
  8 &  8 &   46.193522941 &   37.371345468 &   32.663753877 &   28.549166948\\
  9 & 10 &   52.851064510 &   42.757409793 &   37.371346747 &   32.663753877\\
 10 & 12 &   60.468109856 &   48.919728988 &   42.757411256 &   37.371346747\\
 11 & 15 &   69.182945386 &   55.970179106 &   48.919730662 &   42.757411256\\
 12 & 25 &   79.153787735 &   64.036760096 &   55.970181022 &   48.919730662\\
 13 & 30 &   90.561656168 &   73.265919622 &   64.036762288 &   55.970181022\\
 14 & 40 &  103.613658963 &   83.825211801 &   73.265922130 &   64.036762288\\
 15 & 50 &  118.546753426 &   95.906339123 &   83.825214670 &   73.265922130\\
 16 & 60 &  135.632047826 &  109.728632786 &   95.906342405 &   83.825214670\\
 17 & 80 &  155.179723323 &  125.543034624 &  109.728636542 &   95.906342405\\
 18 &100 &  177.544665266 &  143.636652918 &  125.543038921 &  109.728636542\\
 19 &130 &  203.132906087 &  164.337974808 &  143.636657834 &  125.543038921\\
 20 &150 &  232.408996764 &  188.022829935 &  164.337980433 &  143.636657834\\
 21 &200 &  265.904440681 &  215.121213573 &  188.022836371 &  164.337980433\\
 22 &300 &  304.227342996 &  246.125093134 &  215.121220936 &  188.022836371\\
 23 &500 &  348.073450708 &  281.597339770 &  246.125101559 &  215.121220936\\
 24 &   &  398.238783847 &  322.181947220 &  281.597349409 &  246.125101559\\
 25 &   &  455.634087108 &  368.615723426 &  322.181958248 &  281.597349409\\
 \hline\hline
\end{tabular} 
\end{center}
\end{table}

\begin{table}
\begin{center} 
    \caption
        {New standard sizes for canvas and stretchers rounded up to mm.\label{Table2}}
       \vspace{2mm}        
\begin{tabular}{|c|c|c|ccc|}
\hline\hline 
~~~{No}~~~ &~~{Point}~~ & ~~~{Width}~~~ &\multicolumn{3}{c|}{~~~~~~Height (cm)~~~~~~}\\
&&(cm)&\multicolumn{1}{c}{~~~~~~F~~~~~~}&\multicolumn{1}{c}{~~~~~~P~~~~~~}&\multicolumn{1}{c|}{~~~~~~M~~~~~~}\\
\hline
  1 &  0 &      18.0 &      14.6 &      12.7 &      11.1\\
  2 &  1 &      20.6 &      16.7 &      14.6 &      12.7\\
  3 &  2 &      23.6 &      19.1 &      16.7 &      14.6\\
  4 &  3 &      27.0 &      21.8 &      19.1 &      16.7\\
  5 &  4 &      30.8 &      25.0 &      21.8 &      19.1\\
  6 &  5 &      35.3 &      28.5 &      25.0 &      21.8\\
  7 &  6 &      40.4 &      32.7 &      28.5 &      25.0\\
  8 &  8 &      46.2 &      37.4 &      32.7 &      28.5\\
  9 & 10 &      52.9 &      42.8 &      37.4 &      32.7\\
 10 & 12 &      60.5 &      48.9 &      42.8 &      37.4\\
 11 & 15 &      69.2 &      56.0 &      48.9 &      42.8\\
 12 & 25 &      79.2 &      64.0 &      56.0 &      48.9\\
 13 & 30 &      90.6 &      73.3 &      64.0 &      56.0\\
 14 & 40 &     103.6 &      83.8 &      73.3 &      64.0\\
 15 & 50 &     118.5 &      95.9 &      83.8 &      73.3\\
 16 & 60 &     135.6 &     109.7 &      95.9 &      83.8\\
 17 & 80 &     155.2 &     125.5 &     109.7 &      95.9\\
 18 &100 &     177.5 &     143.6 &     125.5 &     109.7\\
 19 &130 &     203.1 &     164.3 &     143.6 &     125.5\\
 20 &150 &     232.4 &     188.0 &     164.3 &     143.6\\
 21 &200 &     265.9 &     215.1 &     188.0 &     164.3\\
 22 &300 &     304.2 &     246.1 &     215.1 &     188.0\\
 23 &500 &     348.1 &     281.6 &     246.1 &     215.1\\
 24 &   &     398.2 &     322.2 &     281.6 &     246.1\\
 25 &   &     455.6 &     368.6 &     322.2 &     281.6\\

 \hline\hline
\end{tabular} 
\end{center}
\end{table}

\begin{figure}
       \includegraphics[width=17cm]{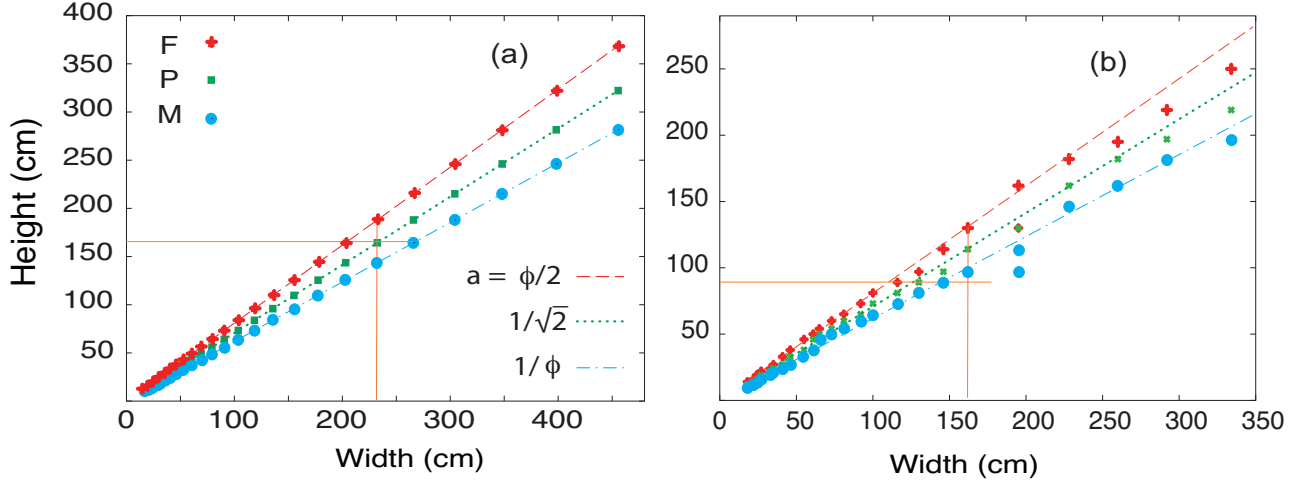}
       \caption{(Color online) Heights of F (red crosses), P (green boxes), and M (blue circles) formats as functions of their widths from the new standard sizes (a) and the French ones including their Japanese extension (b). The (red) dashed, (green) dotted, and (blue) dash-dotted straight lines represent function $y=ax$ with $a=\varphi/2$, $1/\sqrt{2}$, and $1/\varphi$ for the F, P, and M formats, respectively. \label{wh}}
    \end{figure}

Show in Fig. \ref{wh} are the heights $F(i)$, $P(i)$, and $M(i)$  of the F, P and M formats, respectively, plotted as functions of their width $W(i)$ from the new standard sizes [Fig. \ref{wh} (a)] and the French ones [Fig. \ref{wh} (b)]. The dependences of the heights on the width in the new sizes follow the straight line $y=ax$ with $a=\varphi/2$ for the F formats, $a=1/\sqrt{2}$ for the P formats, and $a=1/\varphi$ for the M formats [Fig. \ref{wh} (a)]. The line parallel to the $y$ axis demonstrates that 3 heights $F(i)$, $P(i)$ and $M(i)$ have the same width $W(i)$, whereas the line parallel to the $x$ axis shows that 3 widths $W(i)$, $W(i+1)$, and $W(i+2)$ have the same height, fulfilling the rule for optimal use of material. The French sizes and their Japanese extension do not follow these rules as show in Fig. \ref{wh} (b), where the symbols are not aligned to any straight lines.

As the heights in the new standard sizes are equal to the width $W(i)$ divided by $2/\varphi$ for the F format, $\sqrt{2}$ for the P format, and $\varphi$ for the M format, one immediately sees that the ratio of the canvas circumference $C(i+1)$ to $C(i)$ is equal that of their widths, that is
\begin{equation}
\frac{C(i+1)}{C(i)} = \frac{W(i+1)}{W(i)} = \frac{\varphi}{\sqrt{2}} \simeq 1.144~,
\label{Cratio}
\end{equation}
according to Eq. (\ref{Wi}).

In the same way, one finds the ratio of the canvas area $S(i+1)$ to $S(i)$ to be 
\begin{equation}
\frac{S(i+1)}{S(i)} = \bigg[\frac{W(i+1)}{W(i)}\bigg]^2 = \frac{\varphi^2}{2} \simeq 1.31~.
\label{Sratio}
\end{equation}

From Eqs. (\ref{Cratio}) and (\ref{Sratio}) it follows that
\begin{equation}
C_k(i) = \bigg(\frac{\varphi}{\sqrt{2}}\bigg)^{i-1}C_k(1)\simeq (1.144)^{i-1}C_k(1)~,
\label{Ck}
\end{equation}
\begin{equation}
S_k(i)= \bigg(\frac{\varphi^2}{2}\bigg)^{i-1}S_k(1) \simeq (1.31)^{i-1}S_k(1)~,
\label{Sk}
\end{equation}
with $k= F$, $P$, and $M$, where
\begin{equation}
C_F(1) = 2w(1+\frac{\varphi}{2})\simeq 65.1 {\rm}~{\rm cm}~,\hspace{5mm} S_F(1) = w^2\frac{\varphi}{2}\simeq 262.1~{\rm cm}^2~,
\label{CF1SF1}
\end{equation}
\begin{equation}
C_P(1) = 2w(1+\frac{1}{\sqrt{2}})\simeq 61.5 {\rm}~{\rm  cm}~,\hspace{5mm} S_P(1) = \frac{w^2}{\sqrt{2}}\simeq 229.1~{\rm cm}^2~,
\label{CP1SP1}
\end{equation}
\begin{equation}
C_M(1) = 2w\varphi\simeq 58.3 {\rm}~{\rm cm}~,\hspace{5mm} S_M(1) = \frac{w^2}{\varphi}\simeq 200.2~ {\rm cm}^2~.
\label{CM1SM1}
\end{equation}
In Eq. (\ref{CM1SM1}) the property $1+1/\varphi=\varphi$ of the golden number $\varphi$ is employed. 

The circumferences and areas corresponding to the new standard canvas and stretcher sizes in Table \ref{Table2} are calculated and presented in Tables \ref{Table3} and \ref{Table4}, respectively, so that one can readily use to estimate the cost of any stretcher and pre-stretched canvas, which are manufactured by using the new standard sizes.
\begin{table}
\begin{center} 
    \caption
        {Circumferences of new standard canvas and stretchers rounded up to mm. \label{Table3}}
       \vspace{2mm}        
\begin{tabular}{|c|c|ccc|}
\hline\hline 
~~~{No}~~~ &~~{Point}~~ &\multicolumn{3}{c|}{~~~~~~Circumference (cm)~~~~~~}\\
&&\multicolumn{1}{c}{~~~~~~F~~~~~~}&\multicolumn{1}{c}{~~~~~~P~~~~~~}&\multicolumn{1}{c|}{~~~~~~M~~~~~~}\\
\hline
 1 &  0 &      65.1 &      61.5 &      58.2\\
  2 &  1 &      74.5 &      70.3 &      66.6\\
  3 &  2 &      85.2 &      80.4 &      76.2\\
  4 &  3 &      97.5 &      92.0 &      87.2\\
  5 &  4 &     111.6 &     105.3 &      99.8\\
  6 &  5 &     127.7 &     120.5 &     114.2\\
  7 &  6 &     146.1 &     137.8 &     130.7\\
  8 &  8 &     167.1 &     157.7 &     149.5\\
  9 & 10 &     191.2 &     180.4 &     171.0\\
 10 & 12 &     218.8 &     206.5 &     195.7\\
 11 & 15 &     250.3 &     236.2 &     223.9\\
 12 & 25 &     286.4 &     270.2 &     256.1\\
 13 & 30 &     327.7 &     309.2 &     293.1\\
 14 & 40 &     374.9 &     353.8 &     335.3\\
 15 & 50 &     428.9 &     404.7 &     383.6\\
 16 & 60 &     490.7 &     463.1 &     438.9\\
 17 & 80 &     561.4 &     529.8 &     502.2\\
 18 &100 &     642.4 &     606.2 &     574.5\\
 19 &130 &     734.9 &     693.5 &     657.4\\
 20 &150 &     840.9 &     793.5 &     752.1\\
 21 &200 &     962.1 &     907.9 &     860.5\\
 22 &300 &    1100.7 &    1038.7 &     984.5\\
 23 &500 &    1259.3 &    1188.4 &    1126.4\\
 24 &  &    1440.8 &    1359.7 &    1288.7\\
 25 &   &    1648.5 &    1555.6 &    1474.5\\
 \hline\hline
\end{tabular} 
\end{center}
\end{table}
\begin{table}
\begin{center} 
    \caption
        {Areas of new standard canvas and stretchers rounded up to tenths of cm$^2$. \label{Table4}}
       \vspace{2mm}        
\begin{tabular}{|c|c|ccc|}
\hline\hline 
~~~{No}~~~ &~~{Point}~~ &\multicolumn{3}{c|}{~~~~~~Area (cm$^2$)~~~~~~}\\
&&\multicolumn{1}{c}{~~~~~~F~~~~~~}&\multicolumn{1}{c}{~~~~~~P~~~~~~}&\multicolumn{1}{c|}{~~~~~~M~~~~~~}\\
\hline
   1 &  0 &     262.1 &     229.1 &     200.2\\
  2 &  1 &     343.1 &     299.9 &     262.1\\
  3 &  2 &     449.2 &     392.6 &     343.1\\
  4 &  3 &     587.9 &     513.9 &     449.2\\
  5 &  4 &     769.6 &     672.7 &     587.9\\
  6 &  5 &    1007.5 &     880.6 &     769.6\\
  7 &  6 &    1318.8 &    1152.7 &    1007.5\\
  8 &  8 &    1726.3 &    1508.9 &    1318.8\\
  9 & 10 &    2259.8 &    1975.1 &    1726.3\\
 10 & 12 &    2958.1 &    2585.5 &    2259.8\\
 11 & 15 &    3872.2 &    3384.4 &    2958.1\\
 12 & 25 &    5068.8 &    4430.3 &    3872.2\\
 13 & 30 &    6635.1 &    5799.3 &    5068.8\\
 14 & 40 &    8685.4 &    7591.4 &    6635.1\\
 15 & 50 &   11369.4 &    9937.2 &    8685.4\\
 16 & 60 &   14882.7 &   13008.0 &   11369.4\\
 17 & 80 &   19481.7 &   17027.7 &   14882.7\\
 18 &100 &   25501.9 &   22289.5 &   19481.7\\
 19 &130 &   33382.5 &   29177.3 &   25501.9\\
 20 &150 &   43698.2 &   38193.6 &   33382.5\\
 21 &200 &   57201.7 &   49996.1 &   43698.2\\
 22 &300 &   74878.0 &   65445.8 &   57201.7\\
 23 &500 &   98016.6 &   85669.6 &   74878.0\\
 24 &   &  128305.3 &  112143.0 &   98016.6\\
 25 &   &  167953.9 &  146797.1 &  128305.4\\
 \hline\hline
\end{tabular} 
\end{center}
\end{table}

    \begin{figure}
       \includegraphics[width=12cm]{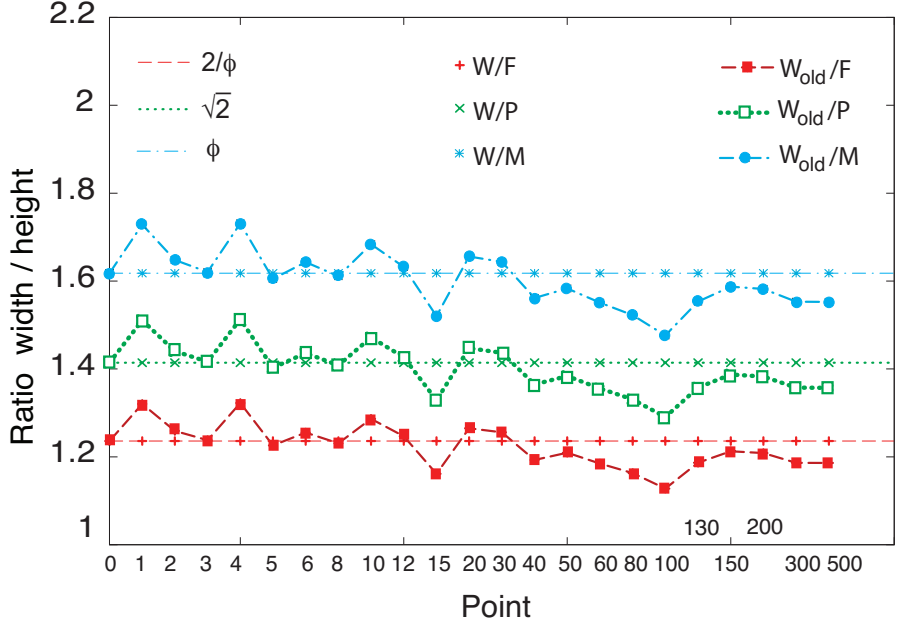}
       \caption{(Color online)  Width-to-height ratios as functions of points.  The (red) crosses, (green) saltires, and (blue) asterisks denote the ratios of the new widths to the new heights in the  F, P, and M formats, respectively. The (red) full boxes, (green) open boxes, and (blue) circles stand for the ratios of the old widths $W_{\rm old}$ to these new heights. The lines connecting the symbols are drawn to guide the eye. The horizontal straight lines stand for $2/\varphi$ (red dashed), $\sqrt{2}$ (green dotted), and $\varphi$ (dash-dotted) as in Fig. \ref{ratioF}. \label{rgoldF}}
    \end{figure}
When the new widths in Tables \ref{Table3} and \ref{Table4} are assigned the points corresponding to the widths with a close length from the French standard sizes, the point 20 and 120 are left out as there are no similar lengths in the new widths. Therefore, in Tables \ref{Table3} and \ref{Table4} there are only 23 sizes corresponding to the French standard widths and their Japanese extension, which are denoted as $W_{\rm old}(i)$  herewith. It is interesting to see how the ratios of the old widths $W_{\rm old}(i)$ to the new F, P and M heights at the same point sizes from Tables \ref{Table3} and \ref{Table4} will deviate from the golden ratio and the porte d'harmonie number. As can be seen from Fig. \ref{rgoldF}, these deviations (from the horizontal straight lines) turn out to be much smaller compared to those of the French standard sizes and their Japanese extension shown in Fig. \ref{ratioF}. From 23 sizes of each format, 12 (that is 52\%) are very close to the corresponding horizontal straight lines. This shows that the new heights, obtained based on the new widths, are even more compatible with the old widths in the sense of approximating the golden ratio and porte d'harmonie number. This is another sound evidence that the French standard canvas sizes and their Japanese extension are not based on the symmetries such as golden and silver ratios as often claimed.
\section{Conclusions}
In the present work a formula has been derived to calculate the widths (longer lengths) of canvases and stretchers, whose heights (shorter lengths) in three figure, landscape and marine formats satisfy the golden ratio and the porte d'harmonie number as well as the rule for optimal use of material.

By using this formula the new standard canvas and stretcher sizes as well as their corresponding circumferences and areas have been calculated starting from the smallest width, which is equal to that of the French standards. 

The exact fulfillment of the golden ratio and the porte d'harmonie number is the main feature that makes the new standard canvas sizes distinctive from the existing French, European and American ones. At the same time, the rule for optimal use of material, which the new sizes also follow, allow them to be easily manufactured in mass production.

Finally, the fact that the new heights are more compatible with the old widths from the French standard canvas sizes and their Japanese extension with respect to the golden ratio and porte d'harmonie number is another sound evidence
that the old standard sizes are not based on these symmetries.

\acknowledgements
The author is thankful to Dr. E. Schwartz of l'\'{E}cole Nationale Sup\'{e}rieure des Beaux-Arts de Paris for his valuable assistance and helpful suggestions.

Thanks are also due to Prof. A. Callen (Australian National University at Canberra) and Dr. P. Labreuche (Paris) for valuable exchanges of information, useful comments, and fruitful discussions.


\begin{thebibliography}{99}
\bibitem{England}J. Simon, {\it Three-quarters, kit-cats and half-lengths: British portrait painters and their canvas sizes, 1625-1850}, National Portrait Gallery, April 2013,\\
http://www.npg.org.uk/research/programmes/artists-their-materials-and-suppliers/three-quarters-kit-cats-and-half-lengths-british-portrait-painters-and-their-canvas-sizes-1625-1850.php
\bibitem{Labreuche1}P. Labreuche, {\it The industrialization of artists' prepared canvas in nineteenth century Paris. Canvas and stretchers: Technical developments up to the period of Impressionism}, International Symposium at the Wallraf-Richatz-Museum \& Fondation Corboud Cologne ``Latest research into painting techniques of Impressionists and Post-Impressionists", 12 - 14, 2008, Zeitschrift f\"{u}r Kunsttechnologie und Konservierung, 22. Jahrgang 2008 Helft 2, pp. 316 - 328.
\bibitem{Labreuche2}P. Labreuche, {\it Paris, capitale de la toile \`{a} peindre, XVIIIe-XIXe si\`{e}cle} (Paris, CTHS/INHA, 2011) pp. 33 - 51.
\bibitem{dynamic}The symmetries based on the irrational numbers such as golden and silver ratios or porte d'harmonie number were classified as dynamic symmetry by American artist Jey Hambridge in his books  ``{\it Dynamic Symmetry: The Greek Vase}" (1920) and ``{\it The Elements of Dynamic Symmetry}" (1926) to distinguish them from static symmetry, which is based on ratios of integers. The present work refrains from using these terms to avoid confusion with the concept of dynamic (dynamical) symmetry of a physical system, whose Hamiltonian can be expressed in terms of Casimirs operators of a chain of algebras.
\bibitem{MAU} Standard Canvas Sizes, Musashi Art \& Design glossary, Musashi Art \& Design University,
http://art-design-glossary.musabi.ac.jp/standard-canvas-sizes/.
\bibitem{Holbein}The Holbein Technical information Q\&A page at\\ http://www.holbein-works.co.jp/technical\_info/2008/10/fmp.html 
\bibitem{Callen1} A. Callen, PhD thesis (Courtauld Institute of Art, University of London, 1980), Chapter 2, pp. 50-77 (unpublished), and private communication.
\bibitem{Callen2} A. Callen, {\it The art of impressionism: Painting technique and the making of modernity} (Yale University Press, 2000), pp. 15, 18 - 21, and 221.\bibitem{Ann}A. Hoenigswald, {\it PSG Stretchers and Strainers III. Materials and Equipment: 3. Standard Stretcher Sizes}, April 2007,http://www.conservation-wiki.com/wiki/\\PSG\_Stretchers\_and\_Strainers\_-\_III.\_Materials\_and\_Equipment
\bibitem{Pernety}A.J. Pernety, {\it Dictionnaire portatif de peinture, sculpture et gravure avec un trait\'{e} pratique des diff\'{e}rentes manieres de peindre, dont la th\'{e}ories est d\'{e}velopp\'{e}e dans les articles qui en sont susceptibles. Ouvrage utile aux artistes, aux \'{e}l\`{e}ves et aux amateurs} (Paris, Bauche Librairie, 1757).
\bibitem{market}See, for instance\\
 http://www.artsupply.com/Golden-Ratio-Canvas-from-Masterpiece\_c\_3017.html
\end{thebibliography}
\end{document}